\documentstyle [12pt,epsf, graphicx, amssymb]{article}

\begin{document}
\def\l{\lambda}
\def\m{\mu}
\def\a{\alpha}
\def\b{\beta}
\def\g{\gamma}
\def\d{\delta}
\def\e{\epsilon}
\def\o{\omega}
\def\O{\Omega}
\def\v{\varphi}
\def\t{\theta}
\def\r{\rho}
\def\bs{$\blacksquare$}
\def\bp{\begin{proposition}}
\def\ep{\end{proposition}}
\def\bt{\begin{th}}
\def\et{\end{th}}
\def\be{\begin{equation}}
\def\ee{\end{equation}}
\def\bl{\begin{lemma}}
\def\el{\end{lemma}}
\def\bc{\begin{corollary}}
\def\ec{\end{corollary}}
\def\pr{\noindent{\bf Proof: }}
\def\note{\noindent{\bf Note. }}
\def\bd{\begin{definition}}
\def\ed{\end{definition}}
\def\C{{\mathbb C}}
\def\P{{\mathbb P}}
\def\Z{{\mathbb Z}}
\def\d{{\rm d}}
\def\deg{{\rm deg\,}}
\def\deg{{\rm deg\,}}
\def\arg{{\rm arg\,}}
\def\min{{\rm min\,}}
\def\max{{\rm max\,}}

% MATH -----------------------------------------------------------
\newcommand{\norm}[1]{\left\Vert#1\right\Vert}
\newcommand{\abs}[1]{\left\vert#1\right\vert}

\newcommand{\set}[1]{\left\{#1\right\}}
\newcommand{\setb}[2]{ \left\{#1 \ \Big| \ #2 \right\} }

\newcommand{\IP}[1]{\left<#1\right>}
\newcommand{\Bracket}[1]{\left[#1\right]}
\newcommand{\Soger}[1]{\left(#1\right)}

\newcommand{\Integer}{\mathbb{Z}}
\newcommand{\Rational}{\mathbb{Q}}
\newcommand{\Real}{\mathbb{R}}
\newcommand{\Complex}{\mathbb{C}}

\newcommand{\eps}{\varepsilon}
\newcommand{\To}{\longrightarrow}
\newcommand{\varchi}{\raisebox{2pt}{$\chi$}}

\newcommand{\E}{\mathbf{E}}
\newcommand{\Var}{\mathrm{var}}

% QED box --------------------------------------------------------
\def\squareforqed{\hbox{\rlap{$\sqcap$}$\sqcup$}}
\def\qed{\ifmmode\squareforqed\else{\unskip\nobreak\hfil
\penalty50\hskip1em\null\nobreak\hfil\squareforqed
\parfillskip=0pt\finalhyphendemerits=0\endgraf}\fi}

% This Document only ---------------------------------------------
\renewcommand{\th}{^{\mathrm{th}}}
\newcommand{\Dif}{\mathrm{D_{if}}}
\newcommand{\Difp}{\mathrm{D^p_{if}}}
\newcommand{\GHF}{\mathrm{G_{HF}}}
\newcommand{\GHFP}{\mathrm{G^p_{HF}}}
\newcommand{\f}{\mathrm{f}}%fitting data
\newcommand{\fgh}{\mathrm{f_{gh}}}%fitting polynomial
\newcommand{\T}{\mathrm{T}}%Taylor polynomial vector
\newcommand{\K}{^\mathrm{K}}%Orde of the Taylor polynomial
\newcommand{\PghK}{\mathrm{P^K_{f_{gh}}}}%the generalized hermite interpolating polynomioal
\newcommand{\Dig}{\mathrm{D_{ig}}}%the digonal matrix with the singular values on the diagonal
\newcommand{\for}{\mathrm{for}}
\newcommand{\End}{\mathrm{end}}

\newtheorem{th}{Theorem}[section]
\newtheorem{lemma}{Lemma}[section]
\newtheorem{definition}{Definition}[section]
\newtheorem{corollary}{Corollary}[section]
\newtheorem{proposition}{Proposition}[section]

\begin{titlepage}

\begin{center}

\topskip 5mm

{\LARGE{\bf {``Smooth rigidity'' and Remez-type inequalities}}}

\vskip 8mm

{\large {\bf Y. Yomdin}}

\vspace{6 mm}
\end{center}

{Department of Mathematics, The Weizmann Institute of Science,
Rehovot 76100, Israel}

\vspace{6 mm}

{e-mail: yosef.yomdin@weizmann.ac.il}

\vspace{1 mm}

\vspace{6 mm}
\begin{center}

{ \bf Abstract}
\end{center}

{\small If a $(d+1)$-smooth function $f(x)$ on $[-1,1]$, with $\max_{[-1,1]}|f(x)|\ge 1,$ has $d+1$ or more distinct zeroes on $[-1,1]$, then
$\max_{[-1,1]}|f^{(d+1)}(x)|\ge 2^{-d-1}(d+1)!$. This follows from the polynomial interpolation of $f$ at its zeroes, with Lagrange's remainder formula.

This is one of the simplest examples of what we call ``smooth rigidity'': certain geometric properties of zero sets of smooth functions $f$ imply explicit lower bounds on the high-order derivatives of $f$. In dimensions greater than one, the powerful one-dimension tools, like Lagrange's remainder formula, and divided finite differences, are not directly applicable. Still, the result above implies, via line sections, rather strong restrictions on zeroes of smooth functions of several variables (\cite{Yom1}).

\smallskip

In the present paper we study the geometry of zero sets of smooth functions, and significantly extend the results of \cite{Yom1}, including into consideration, in particular, finite zero sets (for which the line sections usually do not work). Our main goal is to develop a truly multi-dimensional approach to smooth rigidity, based on polynomial Remez-type inequalities (which compare the maxima of a polynomial on the unit ball, and on its subset). Very informally, {\it one of our main results is that a ``smooth rigidity'' of a zeroes set $Z$ is approximately the  ``inverse Remez constant'' of $Z$}.}

\vspace{1 mm}
\begin{center}
------------------------------------------------
\vspace{1 mm}
\end{center}
This research was supported by the Minerva Foundation.

\end{titlepage}

\newpage

%%%%%%%%%%%%%%%%%%%%%%%%%%%%%%%%%%%%%%%%%%%%%%%%%%%%%%%%%%%%%%%%%%%%%
\section{Introduction}
\setcounter{equation}{0}

Let $f(x)$ be a smooth function on the unit $n$-dimensional ball $B^n$.
%Our goal is to compare the behavior of $f$ with that of polynomials of degree $d$.
A ``rigidity inequality'' for $f$ is an explicit lower bound for the $(d+1)$-st derivative of $f$, which holds, if $f$ exhibits certain patterns, forbidden for polynomials of degree $d$.

\medskip

We expect rigidity inequalities to be valid for those polynomial behavior patterns, which are stable with respect to smooth approximations. At present many such important patterns are known. For those directly relevant to the present paper see \cite{Yom1,Yom.Com,Yom2}. However, translation of the known ``near-polynomiality'' results into ``rigidity inequalities'' usually is not straightforward, and many new questions arise.
%Conversely, if the $(d+1)$-st derivative of $f$ is bounded from above by a certain specific constant $C>0$, its behavior is ``polynomial-like''.

\medskip

The following example illustrates the patterns we are working with in the present paper. Start with a basic property of polynomials: a nonzero univariate polynomial $P(x)$ of degree $d$ can have at most $d$ real zeros. Here is one of the corresponding ``rigidity'' results (well-known in various forms):

\bp\label{prop:d.zeroes}
For each $(d+1)$-smooth function $f(x)$ on $[-1,1]$, with $\max_{[-1,1]}|f(x)|\ge 1$ and with $d+1$ or more distinct zeroes on $[-1,1]$, we have
$$
\max_{[-1,1]}|f^{(d+1)}(x)|\ge \frac{(d+1)!}{2^{d+1}}.
$$
\ep
\pr Let $f(x)$ be a $(d+1)$-smooth function on $[-1,1]$, vanishing at the points $x_1,\ldots,x_{d+1}$. Consider the interpolating polynomial $P$ of $f$, of degree $d$, on the points $x_1,\ldots,x_{d+1}$. We conclude that $P(x)\equiv 0$. From the polynomial interpolation formula with the remainder $R(x)=f(x)-P(x)=f(x)$, given in the Lagrange form, we have:
$$
1\le \max_{[-1,1]}|f(x)| \le \frac{2^{d+1}}{(d+1)!} \max_{[-1,1]}|f^{(d+1)}(x)|,
$$
or
$$
\max_{[-1,1]}|f^{(d+1)}(x)| \ge  \frac{(d+1)!}{2^{d+1}}.
$$
%We immediately conclude also that any $(d+1)$-smooth function $f(x)$ with $\max_{[-1,1]}|f^(x)|\ge 1$ and $\max_{[-1,1]}|f^{(d+1)}(x)|< \frac{(d+1)!}{2^{d+1}},$ has at most $d$ zeroes in $[-1,1]$.
This complete the proof of Proposition \ref{prop:d.zeroes}. $\square$

\bc\label{cor:zeroes}
Any $(d+1)$-smooth function $f(x)$ on $[-1,1]$, with $\max_{[-1,1]}|f(x)|\ge 1$ and $\max_{[-1,1]}|f^{(d+1)}(x)|< \frac{(d+1)!}{2^{d+1}},$ has at most $d$ zeroes in $[-1,1]$.
\ec

\medskip

In higher dimensions the powerful one-dimension tools: Lagrange's remainder formula, and divided finite differences, are not directly applicable. Still, Proposition \ref{prop:d.zeroes} implies, via line sections, rather strong restrictions on the geometric structure of smooth functions of several variables. In particular, the following result was obtained in \cite{Yom1}:

\bt\label{thm:poly.like.old}(\cite{Yom1})
Let $f(x)$ be a smooth function on the unit $n$-dimensional ball $B^n$, with the $C^0$-norm, equal to one. If for some $d \ge 1$, the norm of the $d+1$-st derivative of $f$ is bounded on $B^n$ by $2^{-d-1}$, then the set of zeroes $Y$ of $f$ is ``similar'' to that of a polynomial of degree $d$.

In particular, $Y$ is contained in a countable union of smooth hypersurfaces, ``many'' straight lines cross $Y$ in not more than $d$ points, and the $n - 1$-volume of $Y$ is bounded by a constant, depending only on $n$ and $d$.
\et
We immediately obtain a corresponding ``rigidity inequality'': if a set $Z\subset B^n$ violates any of these restrictions, then, for each smooth $f$, vanishing on $Z$, the norm of the $d+1$-st derivative of $f$ on $B^n$ is at least $2^{-d-1}$.

\smallskip

Still, with one-dimensional tools of Theorem \ref{thm:poly.like.old} we cannot provide any information even for $Z$ being a finite set (unless, by a rare coincidence, many points of $Z$ are on the same straight line).

\medskip

In the present paper we concentrate, as above, on the geometry of the zero set of $f$, and significantly extend the results of Theorem  \ref{thm:poly.like.old}, including into consideration, in particular, finite subsets of zero sets.

\smallskip

Our main goal is to develop a truly multi-dimensional approach to smooth rigidity, based on polynomial Remez-type inequalities (which compare the maxima of a polynomial on the unit ball, and on its subset). Very informally, {\it one of our main results is that ``smooth rigidity'' of a zero set $Z$ is approximately the ``inverse Remez constant'' of $Z$}.

\smallskip

In study of the rigidity of $Z$, as well as of its Remez constant, many approaches are possible: one can stress algebraic geometry of $Z$, its topology, or even its arithmetic geometry (see \cite{Bru.Yom}). We stress the metric geometry of $Z$, and mostly one can think below of finite $Z$. This agrees with a powerful approach to the Whitney smooth extension problem (\cite{Whi1}-\cite{Whi3}), developed in \cite{Bru.Shv,Fef,Kla.Fef} and in many other related publications.

\medskip

In Section \ref{Sec:zero.sets} we define the ``smooth rigidity'', and give some initial examples. In Section \ref{Sec:Remez.zero.sets} the same is done for the Remez constant. Section \ref{Sec:main.results} presents the main results, and Section \ref{Sec:proofs} their proofs.

%It is well known that {\it any close subset $Z\subset B^n$} can be a set of zeroes of a certain $C^\infty$-smooth function on $B^n$. Therefore, one can expect here a certain effect of ``phase transition'': as the higher derivatives of $f$ become smaller, the possible shape of the zero set of $f$ changes from ``any closed'' sets to ``algebraic-like'' sets, as above. In section \ref{Sec:zero.sets} below we provide some results and questions in this direction.

\section{$d$-Rigidity of a zero set $Z$}\label{Sec:zero.sets}
\setcounter{equation}{0}

Let $f: B^n \rightarrow {\mathbb R}$ be a $d+1$ times continuously differentiable function on $B^n$. For $l=0,1,\ldots,d+1$ put
$$
M_l(f)=\max_{z\in B^n} \Vert f^{(l)}(z) \Vert,
$$
where the norm of the $l$-th derivative $f^{(l)}(z)$  of $f$ is defined as the sum of the absolute values of all the partial derivatives of $f$ of order $l$.

\medskip

For $Z\subset B^n$ let $U_d(Z)$ denote the set of $C^{d+1}$ smooth functions $f(z)$ on $B^n$, vanishing on $Z$, with $M_0(f)=1$.

\bd\label{def:rigidity}
For $Z\subset B^n$ we define the $d$-th rigidity constant ${\cal RG}_d(Z)$ as
$$
{\cal RG}_d(Z)=\inf_{f\in U_d(Z)}M_{d+1}(f).
$$
\ed
By this definition we get immediately $M_{d+1}(f)\ge {\cal RG}_d(Z)$ for any $f(z)$ on $B^n$, vanishing on $Z$, with $M_0(f)=1$. Our goal is to estimate ${\cal RG}_d(Z)$ in terms of accessible geometric features of $Z$.

\medskip

Notice that we do not insist on $Z$ being exactly the set of zeroes $Y(f)$ of the functions $f\in U_d(Z)$, but just require $Z\subset Y(f)$.

\medskip

As an example, consider the case of dimension $n=1$. Here the following important fact holds:

\bp\label{prop:d.points}
For any $Z\subset B^1$ we have ${\cal RG}_d(Z)\ge \frac{(d+1)!}{2^{d+1}},$ if $Z$ consists of at least $d+1$ different points, and ${\cal RG}_d(Z)=0$ if $Z$ consists of at most $d$ different points.
\ep
\pr
If $Z$ consists of at most $d$ different points, then there is a polynomial of degree $d$, vanishing on $Z$. Hence ${\cal RG}_d(Z)=0$. If $Z$ consists of at least $d+1$ different points, the result follows from Proposition \ref{prop:d.zeroes}. $\square$

\medskip

Thus in dimension one the minimal non-zero value of ${\cal RG}_d(Z)$ is $\frac{(d+1)!}{2^{d+1}}.$ One of the results of this paper is that this is not true any more in higher dimensions: for $Z\subset B^n, \ n\ge 2,$ the $d$-rigidity ${\cal RG}_d(Z)$ attains arbitrarily small positive values.

\medskip

For any $Z\subset B^n$ its $d$-rigidity ${\cal RG}_d(Z)$ is finite: it is bounded from above by $M_d(f)$ for any $C^\infty$ function $f$, which vanishes exactly on the closure $\bar Z$ of $Z$. The existence of such functions is a classical fact. It is not difficult to give an explicit upper bound for ${\cal RG}_d(Z)$. In dimension one, the complement of any closed set $Z$ is a countable union of open disjoint intervals $U_j$. Let $L$ be the maximal length of $U_j$. Then ${\cal RG}_d(Z)\le C(\frac{1}{L})^{d+1}.$
%Let $\phi(t)$ be a $C^\infty$ function on $[-1,1],$ positive

\medskip

In dimensions $n>1$ we replace $L$ with the maximal size of the cubes in the Whitney covering of the complement of $Z$ (see \cite{Whi2,Whi3}), and get, essentially, the same upper bound.

\medskip

On the other hand, one can easily construct subsets $Z\subset B^n$ with an arbitrarily big $d$-rigidity: consider a subset $Z_h$ of $[-1,1]$ consisting of the grid-points $x_i=-1+ih, \ i=0,1,\ldots, [\frac{2}{h}]$, with $h\ll 1$.

\bp\label{prop:d.points1}
${\cal RG}_d(Z_h)\ge \frac{(d+1)!}{(2(d+2))^{d+1}}(\frac{1}{h})^{d+1}.$
\ep
\pr
Let $f$ be a smooth function on $[-1,1]$, vanishing at the points of $Z_h$, with $M_0(f)=1$. Let $|f(x_0)|=1, \ x_0\in [-1,1]$. Consider a certain subinterval $I\subset [-1,1]$ of the length $(d+2)h$, containing $x_0$. Application of Proposition \ref{prop:d.zeroes}, re-scaled to the interval $I$, completes the proof. $\square$

\medskip

This construction can be easily extended to higher dimensions. Defining the subset $Z^n_h\subset B^n$ as the union of the concentric spheres in $B^n$ with the radii $r_i=ih$, and restricting considerations to the radial straight line, passing through a point $x_0$ with $|f(x_0)|=1$, we obtain the same lower bound for ${\cal RG}_d(Z^n_h)$, as in Proposition \ref{prop:d.points1}.

\medskip

Another simple observation is the following:

\bp\label{prop:Z.interior}
For any $Z \subset B^n$ with a non-empty interior,
$$
{\cal RG}_d(Z) \ge \frac{(d+1)!}{2^{d+1}}.
$$
\ep
\pr
Let $f\in U_d(Z)$. Fix a certain point $x_1$ with $|f(x_1)|=1$, fix $x_2$ in the interior of $Z$, and let $\ell$ be the straight line through $x_1,x_2$. The restriction $\bar f$ of $f$ to $\ell$ (or, more accurately, to the intesection of $\ell$ with $B^n$) has an entire interval of zeroes near $x_2$, and it satisfies $M_0(\bar f)=1$. Applying Proposition \ref{prop:d.zeroes} to $\bar f$ completes the proof. $\square$

\section{Remez constant of $Z$}\label{Sec:Remez.zero.sets}
\setcounter{equation}{0}

Another ingredient we need is a definition and some properties of the Remez (or Lebesgue, or norming, ...) constant.

\bd\label{Remez.constant}
For a set $Z\subset B^n \subset {\mathbb R}^n$ and for each $d\in {\mathbb N}$ the Remez constant ${\cal R}_d(Z)$ is the minimal $K$
for which the inequality
$$
\sup_{B^n}\vert P \vert \leq K \sup_{Z}\vert P \vert
$$
is valid for any real polynomial $P(x)=P(x_1,\dots,x_n)$ of degree $d$.
\ed
%Thus for each polynomial $P$ of degree $d$ we have
%\be\label{eq:remez.ineq}
%\sup_{B^n}\vert P \vert \leq {\cal R}_d(Z) \sup_{Z}\vert P \vert.
%\ee
Clearly, we always have ${\cal R}_d(Z)\ge 1.$ For some $Z$ the Remez constant ${\cal R}_d(Z)$ may be equal to $\infty$. In fact, ${\cal R}_d(Z)$ is infinite if and only if $Z$ is contained in the set of zeroes
$$
Y_P=\{x\in {\mathbb R}^n, \ | \ P(x)=0\}
$$
of a certain polynomial $P$ of degree $d$. Sometimes it is convenient to use the inverse Remez constant $\hat {\cal R}_d(Z):=\frac{1}{{\cal R}_d(Z)}.$

\subsection{``Remez-type'' inequalities}\label{Sec:Remez.type}

``Remez-type'' inequalities provide an upper bound for ${\cal R}_d(Z)$ in terms of various ``computable'' characteristics of $Z$. In particular, the multi-dimensional Remez inequality (\cite{Bru.Gan}, \cite{Rem}, see also \cite{Erd}) uses the Lebesgue measure of $Z$. It reads as follows:

\bt\label{Remez.1}
For any measurable $Z\subset B^n$ we have
\be\label{Remez.ineq.n}
{\cal R}_d(Z) \ \leq \ T_d ({{1+(1-\lambda)^{1\over n}}\over {1-(1-\lambda)^{1\over n}}})\le (\frac{4n}{\lambda})^d.
\ee
Here $T_d(t)=cos(d \ arccos \ t)$ is the Chebyshev polynomial of degree $d$, \ $\lambda= {{m_n(Z)}\over {m_n(B^n)}},$ with $m_n$ being the Lebesgue measure on ${\mathbb R}^n$. This inequality is sharp and for $n=1$ it coincides with the classical Remez inequality of \cite{Rem}.
\et

%Various ``Remez-type'' inequalities provide an upper bound for ${\cal R}_d(Z)$ in terms of various characteristics of $Z$.

Some other examples and a more detailed discussion can be found in \cite{Bru,Bru.Yom, Yom}. In particular, we will use a result of \cite{Yom} (Theorem \ref{Thm:Remez.type} below), showing that the Lebesgue measure can be replaced in Theorem \ref{Remez.1} with a more sensible geometric invariant $\o(Z)$, allowing us to distinguish between discrete and even finite sets of different geometry. Definition of $\o(Z)$ and some of its properties are given in \cite{Yom}, and shortly stated in the next section.

\bt\label{Thm:Remez.type} (\cite{Yom})
For any measurable $Z \subset B^n$, and for $\lambda= {{m_n(Z)}\over {m_n(B^n)}},$ \ $\bar\lambda= {{\o_d(Z)}\over {m_n(B^n)}},$ we have
$$
(\frac{\lambda}{4n})^d \ \le \ (\frac{\bar \lambda}{4n})^d \ \le \ \hat {\cal R}_d(Z).
$$
\et

\subsubsection{Definition and properties of $\omega_d(Z)$}\label{Sec:defin.omega}

To define $\omega_d(Z)$ let us recall that the covering number $M(\e,A)$ of a metric space $A$ is the minimal number of closed
$\e$-balls covering $A$. Below $A$ will be subsets of ${\mathbb R}^n$ equipped with the $l^\infty$ metric. So the $\e$-balls in
this metric are the cubes $Q^n_\e$.

\smallskip

For $\e>0$ we denote by $M_{n,d}(\e)$ (or shortly $M_d(\e)$) the ``Vitushkin polynomials'' of degree $d-1$ in ${1\over \e}$, which appear in bounding the covering number of the sub-level subset of polynomials (see \cite{Fri.Yom,Vit1,Yom,Yom.Com}):

\be\label{Md}
M_d(\e)=\sum_{i=0}^{n-1}C_i(n,d)({1\over \e})^i.
\ee
In particular,
$$
M_{1,d}(\e)=d, \ M_{2,d}(\e)=(2d-1)^2 + 8d({1\over \e}).
$$
Starting with $n=3$ we use the following expression for $M(n,d)$, given in \cite{Fri.Yom}:

\be\label{eq:Vitushkin.coef}
C_i(n,d)=2^{i}(^n_i)(d-i)^i.
\ee

Now for each subset $Z\subset B^n$ (possibly discrete or finite) we introduce the quantity $\omega_d(Z)$ via the following
definition:

\bd\label{omegad}
Let $Z$ be a subset in $B^n\subset {\mathbb R}^n$. Then $\omega_d(Z)$ is defined as
\be\label{omga.eq}
\omega_d(Z) = \sup_{\e>0} {\e}^n[M(\e,Z)- M_d(\e)].
\ee
\ed

The following properties of $\o_d(Z)$ are obtained in \cite{Yom}:

\medskip

1. For a measurable $Z$ \ \ $\o_d(Z)\geq m_n(Z).$

\smallskip

2. For any set $Z\subset B^n$ the quantities $\o_d(Z)$ form a non-increasing sequence in $d$.

\smallskip

3. For a set $Z$ of Hausdorff dimension $n-1$, if the Hausdorff $n-1$ measure of $Z$ is large enough with respect to $d$, then
$\o_d(Z)$ is positive.

\medskip

More information can by found in \cite{Yom,Yom2} (see also \cite{Fav}).
%Some more detailed calculations of $\o_d(Z)$ are given in Section \ref{Sec:examples} below.

\section{Main results}\label{Sec:main.results}
\setcounter{equation}{0}

An important initial observation, connecting the Remez and rigidity constants is

\bl\label{lem:R.is.Inf}
${\cal RG}_d(Z)=0$ if and only if $\hat {\cal R}_d(Z)=0.$
\el

Our main results show a much deeper connection between the rigidity ${\cal RG}_d(Z)$ and the Remez constant $\hat {\cal R}_d(Z)$:

\bt\label{thm:main1}
For any $Z \subset B^n$, \ \ \ $\frac{(d+1)!}{2}\hat {\cal R}_d(Z)\le {\cal RG}_d(Z)$.
\et
This lower bound is valid for any $Z$, and it is sharp, up to constants, for finite sets, as Theorem \ref{thm:main2} below shows. However, we cannot expect {\it an upper bound} of the form

\be\label{eq:both.sides1}
{\cal RG}_d(Z)\le C(n,d) \hat {\cal R}_d(Z)
\ee
to be valid in general: indeed, by Proposition \ref{prop:Z.interior}, for any $Z \subset B^n$ with a non-empty interior, ${\cal RG}_d(Z) \ge \frac{(d+1)!}{2^{d+1}}.$ On the other hand, sets $Z$ with a non-empty interior may have arbitrarily small Remez constant $\hat {\cal R}_d(Z)$.
For example, let $P$ be a polynomial of degree $d$ with $M_0(P)=1$, and let $Z$ be the $\gamma$-sublevel set of $P$: \ $Z=\{z\in B^n, \ |P(z)|\le \gamma \}$. Clearly, we have $\hat {\cal R}_d(Z) \le \gamma$.

\medskip

Still, for some important types of sets $Z$ the bound (\ref{eq:both.sides1}) holds. In this paper we prove it for finite sets $Z$:

\bt\label{thm:main2}
Let $Z \subset B^n$ be a finite set, and let $\rho$ be the minimal distance between the points of $Z$. Then
$$
\frac{(d+1)!}{2}\hat {\cal R}_d(Z)\le {\cal RG}_d(Z)\le \frac{C(n,d)}{\rho^{d+1}}\hat {\cal R}_d(Z).
$$
\et
This theorem can be considered as a generalization of Proposition \ref{prop:d.points} to higher dimensions. However, in dimension one the upper bound of Theorem \ref{thm:main2} requires explanation. Indeed, assume that the cardinality $|Z|$ of $Z$ is $d+1$ or higher. In this case, by Proposition \ref{prop:d.points}, we have ${\cal RG}_d(Z) \ge \frac{(d+1)!}{2^{d+1}}.$ For $\hat {\cal R}_d(Z)=h\ll 1$, if the minimal distance $\rho$ between the points of $Z$ were ``big'', this could contradict the upper bound of Theorem \ref{thm:main2}. However, by the Cartan lemma, for $\hat {\cal R}_d(Z)=h$, we have $\rho \sim h^{\frac{1}{d}}$, and the upper bound of Theorem \ref{thm:main2} is of order $h^{-\frac{1}{d}}\gg 1.$

\smallskip

In dimensions $2$ and higher we have finite sets $Z$ with arbitrarily small $\hat{\cal R}_d(Z),$ and with $\rho$, uniformly bounded from below. For such sets the upper bound of Theorem \ref{thm:main2} is meaningful. One of the simplest examples is a plane triangle $Z_h$, defined as
$$
Z_h=\{(-\frac{1}{2},0),(0,h),(\frac{1}{2},0)\}.
$$
Easy computation shows that $\hat {\cal R}_1(Z_h)=\frac{h}{2}$.

\medskip

In our next result we use Remez-type inequalities: Theorem \ref{thm:main1}, combined with Theorem \ref{Thm:Remez.type} gives the following lower bound for the rigidity:

\bc\label{cor:main5}
For any measurable $Z \subset B^n$,
$$
\frac{(d+1)!}{2}(\frac{\lambda}{4n})^d \ \le \ \frac{(d+1)!}{2}(\frac{\bar \lambda}{4n})^d \ \le {\cal RG}_d(Z).
$$
\ec

For finite sets $Z$, only the part of Corollary \ref{cor:main5}, involving $\bar\lambda$, remains relevant. One of the simplest examples is the following:

%In Section \ref{Sec:examples} we give some examples where $\o_d(Z)$, and hence $\bar\lambda(Z),$ can be estimated explicitly. Among the sets for which $\o_d(Z)$, and hence $\bar\lambda(Z),$ can be estimated explicitly, are ``sufficiently dense'' nets in hypersurfaces of a ``sufficiently big'' $n-1$-volume, and ``sufficiently dense'' nets in sets of a positive $n$-measure.
%We provide more details in Section \ref{Sec:examples} below.

%A corollary of the results of \cite{Yom} is the following estimate of $V(Z)$ for finite $Z$ in terms of the simplest geometric parameters of $Z$:

\bt\label{thm:main3}
Let $Z \subset B^n$ be a finite set, and let $\rho$ be the minimal distance between the points of $Z$. Assume that the cardinality $M=|Z|$ satisfies
$M > (4d)^n(\frac{1}{\rho})^{n-1}$. Then
$$
0 < \frac{(d+1)!}{2}\left ( \frac{M\rho^n- (4d)^n\rho}{4n}\right )^d \le {\cal RG}_d(Z).
$$
\et
%This result follows directly from Theorem \ref{thm:main1}, combined with Proposition \ref{prop:main5} of Section \ref{Sec:examples}.
In particular, consider a regular grid $Z$ with the step $h$ in the unit cube $[0,1]^n$ in ${\mathbb R}^n$. This grid contains $M=(\frac{1}{h})^n$ points, and the minimal distance between these points is $h$. Thus for
$$
(\frac{1}{h})^n > (4d)^n(\frac{1}{h})^{n-1}, \ \ or \ \ h < (\frac{1}{4d})^n,
$$
the conditions of Theorem \ref{thm:main3} are satisfied. Taking, say, twice a smaller $h$, we can shift each point of the grid $Z$, in an arbitrary direction, to any distance, not exceeding, say, $\frac{h}{10}$. We get a ``near grid'' $Z'$. The same calculation as above shows that the conditions of Theorem \ref{thm:main3} are still satisfied, while (generically) no straight line crosses $Z'$ at more than two point.

\medskip

\subsection{Some open questions}\label{Sec:questions}

The results above show that the rigidity ${\cal RG}_d(Z)$ and the inverse Remez constant $\hat {\cal R}_d(Z)$ behave in a similar way. In particular, for finite sets with minimal separation $\rho$ between the points uniformly bounded from below, ${\cal RG}_d(Z)$ and $\hat {\cal R}_d(Z)$ are equivalent, up to constants.

\smallskip

Still, an important gap between these two quantities exists: for $Z$ with a non-empty interior, the rigidity is uniformly bounded from below by $\frac{(d+1)!}{2^{d+1}}$, while the inverse Remez constant can be arbitrarily small.

\medskip

{\it Are there other natural classes of sets $Z$, for which such a gap exists, besides the sets with a non-empty interior?

\smallskip

%Is the existence of the gap related to the Hausdorff (Minkivski) dimension of $Z$?

In particular, are there finite sets $Z$ with $\hat {\cal R}_d(Z)$ arbitrarily small, but with ${\cal RG}_d(Z)$ bounded from below by a positive constant?

\smallskip

On the other side, are there other natural classes of sets $Z$, besides the finite ones, with the minimal distance between the point uniformly bounded from below, for which ${\cal RG}_d(Z)$ and $\hat {\cal R}_d(Z)$ are equivalent, up to constants?}

\medskip

The last question is closely related to the Whitney smooth extension problem for the restrictions to $Z$ of polynomials of degree $d$ - compare the proof of Theorem \ref{thm:main2} in the next section.

\section{Proof of main results}\label{Sec:proofs}
\setcounter{equation}{0}

\noindent {\large {\bf  Proof of Lemma \ref{lem:R.is.Inf}}}.

\medskip

We prove a slightly more detailed statement:

\bp\label{prop:rg.vs.r1}
For a subset $Z\subset B^n$, and for each $d=1,2,\ldots$ the rigidity constant ${\cal RG}_d(Z)$ is zero if and only of the Remez constant    ${\cal R}_d(Z)$ is infinity. Both these conditions are equivalent to $Z$ being contained in a zero set of a certain polynomial $P(x)$ of degree $d$.
\ep
\pr
In one direction this is immediate: if $Z$ is contained in a zero set of a certain polynomial $P(x)$ of degree $d$, then $P$ itself, being properly normalized, belongs to $U_d(Z)$, and its $d+1$-st derivative is zero. Hence ${\cal RG}_d(Z)$ is zero. Since $P$ is zero on $Z$ and non-zero identically, we see also that ${\cal R}_d(Z)$ is infinity.

\medskip

Now assume that $Z$ is not contained in a zero set of any polynomial $P(x)$ of degree $d$. Then one can easily show that ${\cal R}_d(Z)$ is finite (see, e.g. \cite{Bru.Yom}).

To see that in this case ${\cal RG}_d(Z)$ is nonzero, assume the contrary: ${\cal RG}_d(Z)=0$. This means that there is a sequence of functions $f_m\in U_d(Z)$, such that $M_{d+1}(f_m)\to 0$ as $m\to \infty$. Applying to $f_m$ the Taylor formula of order $d$ at the origin, we obtain a sequence of polynomials $Q_m$ of degree $d$, with $M_0(f_m-Q_m)\to 0$. This implies that $Q_m$ converge to zero on $Z$, while $M_0(Q_m)\to 1$. But since ${\cal R}_d(Z)$ is finite, this is impossible. The contradiction completes the proof of Lemma \ref{lem:R.is.Inf}. $\square$

\medskip
\medskip
\medskip
\medskip
\medskip

\noindent {\large {\bf Proof of Theorem \ref{thm:main1}}}.

\medskip

We obtain Theorem \ref{thm:main1} as a corollary of a ``Remez-type inequality for smooth functions'', obtained in \cite{Yom2}.

It is clear that Remez inequality of Theorem \ref{Remez.1} cannot be verbally extended to smooth functions: such function $f$ may be identically zero on any given closed set $Z$ (in particular, of a positive measure), and non-zero elsewhere. It was shown in \cite{Yom2} that adding a ``remainder term'' (expressible through the bounds on the derivatives of $f$) provides a generalization of the Remez inequality to smooth functions. Here is the result of \cite{Yom2}, which we need:

\medskip

Let $R_k(f)={1\over {(k+1)!}}M_{k+1}(f)$ be the Taylor remainder term of $f$ of degree $k$ on the unit ball $B^n$.

\bt\label{cor:Taylor1}(\cite{Yom2})
Let $f: B^n \rightarrow {\mathbb R}$ be a $d+1$ times continuously differentiable function on $B^n$, and let a subset $Z\subset B^n$ be given. Then for $L=\max_{x\in Z}|f(x)|$ we have
\be\label{Smooth.Remez.ineq1}
M_0(f)=\max_{x\in B^n} \vert f(x) \vert \leq \min_{k=0,1,\ldots,d} [{\cal R}_k(Z)(L+R_k(f))+R_k(f)].
\ee
%where $R_d(f)={1\over {(d+1)!}}M_{d+1}(f)$ is the Taylor remainder term of $f$ of degree $d$ on the unit ball $B^n$.
\et

On this base we produce in \cite{Yom2} also a corresponding general rigidity result:

\bt\label{Remez.Whitney}(\cite{Yom2})
Let a set $Z\subset B^n$ be given, let $f$ be a $C^k$ function on $B^n$, equal to zero on $Z$ and satisfying $M_0(f)=1$. Then for each $d=0,\ldots,k-1$ we have
$$
M_{d+1}(f)\geq {{(d+1)!}\over {{\cal R}_d(Z)+1}}.
$$
\et

In order to prove Theorem \ref{thm:main1}, it remains to notice that, since always ${\cal R}_d(Z)\ge 1,$ we have
$$
{\cal RG}_d(Z)=\inf_{f\in U_d(Z)}M_{d+1}(f) \ \geq {{(d+1)!}\over {{\cal R}_d(Z)+1}} \geq {{(d+1)!}\over {2{\cal R}_d(Z)}}\ = \ \frac{(d+1)!}{2}\hat {\cal R}_d(Z).
$$
This completes the proof. $\square$

\medskip
\medskip
\medskip
\medskip
\medskip

\noindent {\large {\bf Proof of Theorem \ref{thm:main2}}}.

\medskip

We fix a $C^\infty$ function $\psi=\psi_n$ on $B^n$, such that $\psi(0)=1$, and $\psi$ vanishes near the boundary of $B_n$. Put $C(n,d)=2^{d+1}M_{d+1}(\psi).$

\smallskip

Now let $Z=\{z_1,\ldots,z_m\}$, and let $P$ be a polynomial with $M_0(P)=1$, and with $\max_{Z}|P|= \hat {\cal R}_d(Z)$. We define $f(x)$ by
$$
f(x)=P(x)-\sum_{j=1}^m P(z_j)\psi(\frac{2}{\rho}(x-z_j)).
$$
Each function $\psi_j=\psi(\frac{2}{\rho}(x-z_j))$ is supported in a ball $B_j$ of a radius $\frac{\rho}{2}$ centered at $z_j$, and
$$
M_{d+1}(\psi_j)=(\frac{2}{\rho})^{d+1}M_{d+1}(\psi)=\frac{C(n,d)}{\rho^{d+1}}.
$$
Since the balls $B_j$ are disjoint by the assumptions on $Z$, we conclude that $f$ is a $C^\infty$ function, and
$$
M_{d+1}(f)\le (\max_j |P(z_j)|)\frac{C(n,d)}{\rho^{d+1}} = \frac{C(n,d)}{\rho^{d+1}} \hat {\cal R}_d(Z).
$$
This completes the proof. $\square$

\medskip

%The result of Theorem \ref{thm:main2} can be extended to wider classes of sets $Z$, in particular, to $Z$ being a finite disjoint union of smooth compact submanifolds. The proof requires some tools of differential topology, and we plan to present it separately.

\medskip
\medskip
\medskip
\medskip
\medskip

\noindent {\large {\bf Proof of Theorem \ref{thm:main3}}}.

\medskip

We use the following rough bound for the Vitushkin polynomials $M_{n,d}(\e)$:

\bp\label{Prop:Rough.Vit}
$$
M_{n,d}(\e)\le (4d)^n(\frac{1}{\e})^{n-1}.
$$
\ep
\pr
It follows from the expression for the coefficients of the Vitushkin polynomial $M(n,d)$, given in \cite{Fri.Yom}:
$$
C_i(n,d)=2^{i}(^n_i)(d-i)^i,
$$
via direct calculations. $\square$

\bp\label{prop:main5}
Let $Z \subset B^n$ be a finite set, and let $\rho$ be the minimal distance between the points of $Z$. Assume that the cardinality $M=|Z|$ satisfies $M > (4d)^n(\frac{1}{\rho})^{n-1}$. Then
$$
\hat {\cal R}_d(Z)\ge \left ( \frac{M\rho^n- (4d)^n\rho}{4n}\right )^d >0.
$$
\ep
\pr
We substitute $\e=\rho$ in the definition of $\o_d(Z)$, and get, by Proposition \ref{Prop:Rough.Vit},
$$
\o_d(Z)\ge \rho^n(M-M(n,d)(\rho))\ge M\rho^n-(4d)^n\rho.
$$
Via Theorem \ref{Thm:Remez.type} this completes the proof of Proposition \ref{prop:main5}, and hence, via Theorem \ref{thm:main1}, also the proof of Theorem \ref{thm:main3}. $\square$

\bibliographystyle{amsplain}

\end{document}